\documentclass[12pt]{amsart}
\usepackage{amssymb,latexsym,eufrak,amsmath,amscd, graphicx}

\theoremstyle{plain}
\newtheorem{theorem}{Theorem}[section]
\newtheorem{proposition}[theorem]{Proposition}

\newtheorem{lemma}[theorem]{Lemma}

\theoremstyle{definition}
\newtheorem{definition}[theorem]{Definition}
\newtheorem{remark}[theorem]{Remark}
\newtheorem{example}[theorem]{Example}

\newtheorem{question}[theorem]{Question}

\newtheorem{problem}[theorem]{Problem}
\newtheorem{conjecture}[theorem]{Conjecture}

\newcommand{\RR}{\mathbb{R}}
\newcommand{\NN}{\mathbb{N}}
\newcommand{\A}{\mathbb{A}}
\newcommand{\ZZ}{\mathbb{Z}}
\newcommand{\CC}{\mathbb{C}}
\newcommand{\QQ}{\mathbb{Q}}

\newcommand{\cO}{\mathcal  {O}}
\newcommand{\cI}{\mathcal  {I}}
\newcommand{\fra}{\frak{a}}
\newcommand{\frb}{\frak{b}}

\newcommand{\frmm}{\frak{m}}

\DeclareMathOperator{\Ker}{ker}

\DeclareMathOperator{\mult}{mult}

\DeclareMathOperator{\lct}{lc}

\DeclareMathOperator{\pt}{c}

\begin{document}

\title[F-thresholds and Bernstein-Sato polynomials]  {F-thresholds
and Bernstein-Sato polynomials} 

\author[Mustata]{Mircea Musta\c{t}\v{a}}
\address{Department of Mathematics, University of
Michigan\\ Ann Arbor, MI 48109, USA}  
\email{mmustata@umich.edu}

\author[Takagi]{Shunsuke Takagi}
\address{Faculty of Mathematics, Kyushu University, 6-10-1
Hakozaki, Higashi-ku, Fukuoka-city, 812-8581, Japan}
\email{stakagi@math.kyushu-u.ac.jp}

\author[Watanabe]{Kei-ichi Watanabe}
\address{Department of Mathematics, College of Humanities and Sciences,
Nihon University, Setagaya-Ku, Tokyo 156-0045, Japan} 
\email{watanabe@math.chs.nihon-u.ac.jp}

\maketitle

\section*{Introduction}

We introduce and study invariants of singularities 
in positive characteristic 
 called F-thresholds. 
They give an analogue of the jumping coefficients
of multiplier ideals in characteristic zero. Unlike these, however,
the F-thresholds are not defined via
resolution of singularities, but via
the action of the Frobenius morphism.

We are especially interested in the connection between the invariants of an
ideal $\fra$ in characteristic zero 
and the invariants of the different reductions mod $p$ of $\fra$.
Our main point is that this 
relation depends on arithmetic properties of $p$.
We present several examples, as well as some questions on this topic.
In a slightly different direction, we describe a new connection
between invariants mod $p$ and the roots of the Bernstein-Sato polynomial.

We will restrict ourselves to the case of an ambient smooth variety,
when our invariants have a down-to-earth
description. Let $(R,\frmm)$ be a regular local ring of characteristic
$p>0$. We want to measure the singularities of a nonzero ideal
$\fra\subseteq\frmm$. For every ideal $J\subseteq\frmm$ containing $\fra$
in its radical, and for every $e\geq 1$, we put
$$\nu^J_{\fra}(p^e):=\max\{r\vert\fra^r\not\subseteq J^{[p^e]}\},$$
where $J^{[p^e]}=(f^{p^e}\vert f\in J)$. 
One can check that the limit
$$\pt^J(\fra):=\lim_{e\to\infty}\frac{\nu^J_{\fra}(p^e)}{p^e}$$
exists and is finite. We call this limit the F-threshold of $\fra$
with respect to $J$. When $J=\frmm$, we simply write $\pt(\fra)$
and $\nu_{\fra}(p^e)$. The invariant $\pt(\fra)$ was introduced
in \cite{TW} under the name of F-pure threshold.

In the first section we define these invariants and give
 their basic properties.
The second section is devoted to the connection with the generalized test
ideals introduced by Hara and Yoshida in \cite{HY}.
More precisely, we show that our invariants are the jumping coefficients
for their test ideals.
As it was shown in \cite{HY} that
the test ideals satisfy similar properties with the multiplier ideals
in characteristic zero, it is not surprising
that the F-thresholds behave in a similar way with
the jumping 
coefficients of the multiplier ideals from \cite{ELSV}. 
Such an analogy was also stressed in \cite{TW},
where it was shown that the smallest F-threshold $\pt(\fra)$ behaves 
in the same way as the smallest jumping coefficient in characteristic zero
(known as the log canonical threshold).

We point out 
that it is not known whether the analogue of two basic properties 
of jumping coefficients of multiplier ideals hold in our setting: 
whether $\pt^J(\fra)$ is 
always a rational number and whether the set of all F-thresholds of $\fra$
is discrete.

There are very interesting questions related to the invariants
attached to different 
reductions mod $p$ of a characteristic zero ideal $\fra$.
We discuss these in \S 3. For simplicity, we assume that
$\fra$ and $J$ are ideals in $\ZZ[X_1,\ldots,X_n]$, contained in 
$(X_1,\ldots,X_n)$ and such that $\fra$ 
is contained in the radical of $J$.
Let us denote by $\fra_p$ and $J_p$ the localizations at $(X_1,\ldots,X_n)$
of the images of $\fra$ and $J$, respectively, in 
${\mathbb F}_p[X_1,\ldots,X_n]$.
We want to compare our invariants mod $p$ (which we write as
$\nu^J_{\fra}(p^e)$ and $\pt^J(\fra_p)$) with the characteristic zero
invariants of $\fra$ (more precisely, with the invariants around the origin
of the image $\fra_{\QQ}$
of $\fra$ in $\QQ[X_1,\ldots,X_n]$).

First, let us denote by $\lct_0(\fra)$ the log canonical threshold
of $\fra_{\QQ}$ around the origin. 
It follows from results of Hara and Watanabe 
(see \cite{HW}) that if $p\gg 0$ then
$\pt(\fra_p)\leq\lct_0(\fra)$ and $\lim_{p\to\infty}
\pt(\fra_p)=\lct_0(\fra)$. Moreover, results of 
Hara and Yoshida from \cite{HY} allow the extension of these
formulas 
to higher jumping numbers
(see Theorems~\ref{thm1'} and \ref{thm2'} below for statements).

It is easy to give examples in which 
$\pt(\fra_p)\neq\lct_0(\fra)$
for infinitely many $p$. On the other hand, one conjectures
that there are infinitely many $p$ with $\pt(\fra_p)
=\lct_0(\fra)$.
We give examples in which more is true:
there is a positive integer $N$
such that for $p\equiv 1$ (mod $N$) we have equality
$\pt(\fra_p)=\lct_0(\fra)$. Moreover, in these examples
one can find rational functions
$R_i\in\QQ(t)$ associated to every $i\in\{1,\ldots, N-1\}$ 
relatively prime to $N$,
such that $\pt(\fra_p)=R_i(p)$ whenever 
$p\gg 0$ satisfies $p\equiv i$
(mod $N$). 
It would be interesting to understand better when such a behavior holds.
As the example of a cone over an elliptic curve 
without complex multiplication shows, this can't
hold in general. On the other hand, motivated by our examples
one can speculate that the following holds: there is always a number field
$K$ such that whenever the prime $p$ is large enough and completely split in $K$,
then $\pt(\fra_p)=\lct_0(\fra)$. 

A surprising fact is that our invariants for $\fra_p$ are related 
to the Bernstein-Sato polynomial $b_{\fra,0}(s)$ of $\fra$. 
More precisely, we show
 that for all $p\gg 0$ and for all $e$, 
we have $b_{\fra,0}(\nu^J_{\fra}(p^e))\equiv 0$ (mod $p$).
We show on some examples in \S 4 how to use this to give roots of
the Bernstein-Sato polynomial (and not just roots mod $p$).

In these examples we will see the following behavior:
given some ideal $J$ containing $\fra$ in its radical,
and $e\geq 1$, we can find $N$
such that for all $i\in\{1,\ldots, N-1\}$ 
relatively prime to $N$ there are polynomials
$P_i\in\QQ[t]$ of degree $e$ satisfying $\nu^J_{\fra}(p^e)=P_i(p)$
for all $p\gg 0$, with $p\equiv i$ (mod $N$). 
The previous observation implies that $b_{\fra,0}(P_i(0))$ is divisible by $p$
for every such $p$. By Dirichlet's Theorem we deduce that $P_i(0)$
is a root of $b_{\fra,0}$.

An interesting question is which roots can be obtained by the above method.
It is shown in \cite{BMS1} that for monomial ideals the functions
$p\to\nu^J_{\fra}(p^e)$ behave as described above, and moreover,
all roots of the Bernstein-Sato polynomial are given 
by this procedure.
On the other hand, Example~\ref{quadric} below
shows that in some cases there are roots which can not
be given by our method.

\section{F-thresholds}

Let $(R,\frmm,k)$ be a regular local ring 
of dimension $n$ and of characteristic $p>0$.
Since $R$ is regular, the Frobenius morphism $F : R\longrightarrow R$, 
$F(x)=x^p$
is flat.

In what follows $q$ denotes a positive power of $p$, and if 
$I=(y_1,\ldots,y_s)$ is an ideal in $R$, then 
$$I^{[q]}:=(y^q\vert y\in I)=(y_1^{q},\ldots,y_s^{q}).$$

We will use below
the fact that as $R$ is regular, every ideal $I$ is equal with its tight closure
(see, for example \cite{HH}). This means that if $u$, $f\in R$
are such that $uf^q\in I^{[q]}$ for all $q\gg 0$, and if $u\neq 0$, then
$f\in I$. This is easy to see: by the flatness of the Frobenius morphism
we have $(I^{[q]}:f^q)=(I:f)^{[q]}$. 
Therefore  $u$ lies in $\bigcap_q(I:f)^{[q]}$, which is zero
$f$ is not in  $I$.

Let $\fra$ be a fixed ideal of $R$, such that $(0)\neq\fra\subseteq\frmm$.
To each ideal $J$ of 
$R$ such that $\fra\subseteq {\rm Rad}(J)\subseteq\frmm$,
we associate a threshold as follows. For every $q$, let
$$\nu_{\fra}^J(q):=\max\{r\in\NN\vert\fra^r\not\subseteq J^{[q]}\}.$$ 
As $\fra\subseteq {\rm Rad}(J)$, this is a nonnegative
integer.

\begin{lemma}\label{lem1}
For every $\fra$, $J$ and $q$ as above, we have
$\nu_{\fra}^J(pq)\geq p\cdot\nu_{\fra}^J(q)$.
\end{lemma}

\begin{proof}
The inequality is a consequence of the fact that if $u\not\in J^{[q]}$,
then $u^p\not\in J^{[pq]}$.
\end{proof}

It follows from the above lemma that 
\begin{equation}\label{eq1}
\lim_{q\to\infty}\frac{\nu_{\fra}^J(q)}{q}=\sup_q\frac{\nu_{\fra}^J(q)}{q}.
\end{equation}
We call this limit the \emph{F-threshold} of the pair $(R,\fra)$ (or simply of
$\fra$) with respect to $J$, and we denote it by $\pt^J(\fra)$.

\begin{remark}\label{finite}
The above limit is finite. In fact, if $\fra$ is generated by $r$ elements,
and if $\fra^N\subseteq J$, then 
$$\fra^{N(r(p^e-1)+1)}\subseteq(\fra^{[p^e]})^{N}
=(\fra^N)^{[p^e]}\subseteq J^{[p^e]}.$$
Therefore $\nu_{\fra}^J(p^e)\leq N(r(p^e-1)+1)-1$. Dividing by $p^e$
and taking the limit gives $\pt^J(\fra)\leq Nr$.

We also have $\pt^J(\fra)>0$. More precisely, as $\fra\neq (0)$,
Krull's Intersection Theorem shows that we can find
$e$ such that $\fra\not\subseteq J^{[e]}$, so $c^J(\fra)\geq 1/p^e$.

We make the convention $c^R(\fra)=0$.
\end{remark}

\begin{example}\label{parameter}
If $J$ is an ideal generated by a regular sequence $y_1,\ldots,y_r$
in $R$, then $\nu_J^J(q)=r(q-1)$ for all $q$. Therefore $\pt^J(J)=r$.
\end{example}

\begin{question}\label{rationality}
Is it true that for all 
nonzero ideals $\fra$ and $J$ with $\fra\subseteq
{\rm Rad}(J)\subseteq\frmm$, the F-threshold $c^J(\fra)$ is a rational 
number ?
\end{question}

\smallskip

\begin{remark}
The F-pure threshold $\pt(\fra)$ was defined in \cite{TW}
(under the assumption that the Frobenius morphism $F$ on $R$ is finite) as
the supremum of those $t\in\QQ_+$ such that 
the pair $(R,\fra^t)$ is F-pure. Under this extra assumption on
$F$, since $R$ is regular, the pair $(R,\fra^t)$  
is F-pure if and only 
if for $q\gg 0$ we have  $\fra^{\lfloor t(q-1)\rfloor}
\not\subseteq \frmm^{[q]}$ (see Lemma~3.9 in \cite{Ta}). 
Here we use the notation $\lfloor \alpha\rfloor$
for the largest integer $\leq\alpha$.
 
The above condition is equivalent with $\nu^{\frmm}_{\fra}(q)\geq 
\lfloor t(q-1)\rfloor$ for $q\gg 0$. It follows from our definition
that if $(R,\fra^t)$ is F-pure, then $t\leq \pt^{\frmm}(\fra)$, and that if
$t<\pt^{\frmm}(\fra)$, then $(R,\fra^t)$ is F-pure. Therefore the
F-pure threshold $\pt(\fra)$ is equal to the F-threshold
$\pt^{\frmm}(\fra)$ of $\fra$ with respect to the maximal ideal.
We will keep the notation $\pt(\fra)$ for $\pt^{\frmm}(\fra)$,
and moreover, we will put $\nu_{\fra}(q):=\nu^{\frmm}_{\fra}(q)$.

Note that the F-pure threshold was defined 
in \cite{TW} without the regularity assumption on 
$R$, but in what follows we will work under this restrictive hypothesis.
\end{remark}

\begin{remark}
In characteristic zero, the only analogue of $J^{[q]}$
which does not depend on the choice of generators for $J$ is
the usual power $J^q$. If we imitate the definition of the
F-pure threshold in this setting, replacing $\frmm^{[q]}$ by $\frmm^q$,
then we get $1/\mult_0(\fra)$, where $\mult_0(\fra)$ is the
largest power of $\frmm$ containing $\fra$.
\end{remark}

Here are a few properties of F-thresholds. When $J=\frmm$
these have been proved in \cite{TW} in a more general setting.

\begin{proposition}\label{prop1}

Let $\fra$, $\frb$, $J\subseteq\frmm$ be nonzero ideals,
such that $\fra$ and $\frb$ are contained in the radical of $J$.

\begin{enumerate}
\item If $\fra\subseteq\frb$, then $\pt^J(\fra)\leq\pt^J(\frb)$.
\item $\pt^J(\fra^s)=\frac{\pt^J(\fra)}{s}$ for every positive integer $s$.
\item If $\fra\subseteq J^s$ and $J$ can be generated by $m$ elements,
then $\pt^J(\fra)\leq m/s$. 
If $\fra\not\subseteq \frmm^{s+1}$ and i $J\subseteq
\frmm^{\ell}$, then $\pt^J(\fra)\geq\ell/s$.
\item If $\overline{\fra}$ is the integral closure of $\fra$,
then $\pt^J(\fra)=\pt^J(\overline{\fra})$. 
\item For every $q$, we have $\frac{\nu_{\fra}(q)}{q}<\pt^J(\fra)$.
\item We have $\pt^J(\fra+\frb)\leq\pt^J(\fra)+\pt^J(\frb)$.
\end{enumerate}
\end{proposition}

\begin{proof}
The first assertion is trivial: since $\fra\subseteq\frb$, we get
$\nu^J_{\fra}(q)\leq\nu^J_{\frb}(q)$ for all $q$. Hence $\pt^J(\fra)
\leq\pt^J(\frb)$. 

Given $s$ and $q$, we have $(\fra^s)^r\not\subseteq J^{[q]}$
if and only if $rs\leq\nu^J_{\fra}(q)$. Hence $\nu^J_{\fra^s}(q)=
\lfloor \nu_{\fra}^J(q)/s\rfloor$, 
which after dividing by $q$ and passing to limit gives ${\rm (}2{\rm )}$.

If $\fra\subseteq J^s$, then by ${\rm (}1{\rm )}$ and 
${\rm (}2{\rm )}$ we have 
$\pt^J(\fra)\leq\pt^J(J^s)=\frac{\pt^J(J)}{s}\leq\frac{m}{s}$. 
The last inequality follows from Remark~\ref{finite}.

Suppose now that $\fra\not\subseteq\frmm^{s+1}$. If $\pt^J(\fra)<\ell/s$, 
then by taking
$q$ large enough we can find $r$ such that 
\begin{equation}\label{eq3}
\pt^J(\fra)<\frac{r}{q}<\frac{\ell}{s}.
\end{equation}
The first inequality 
shows that $\fra^r\subseteq J^{[q]}\subseteq\frmm^{q\ell}$.
As by hypothesis $\fra^r\not\subseteq\frmm^{rs+1}$, we deduce
$rs+1>q\ell$. This contradicts 
the second inequality in (\ref{eq3}).

For ${\rm (}4{\rm )}$ note first that $\pt^J(\fra)\leq\pt^J(\overline{\fra})$
follows from ${\rm (}1{\rm )}$. For the reverse inequality, 
recall that by general properties of the integral closure,
there is a fixed positive integer $s$ such that
$\overline{\fra}^{\ell+s}
\subseteq\fra^{\ell}$ for all $\ell$. 
Hence we have $\nu^J_{\fra}(q)\geq \nu^J_{\overline{\fra}}(q)-s$
for every $q$, which implies $\pt^J(\fra)\geq\pt^J(\overline{\fra})$.

In order to prove ${\rm (}5{\rm )}$ suppose that for some $q$ we have
$\nu^J_{\fra}(q)/q=\pt^J(\fra)$. If $\nu^J_{\fra}(q)=r$, this implies
that $\nu^J_{\fra}(qq')=rq'$ for all $q'$.
Therefore $\fra^{rq'+1}\subseteq J^{[qq']}$ for all $q'$.
As $J^{[q]}$ is equal to its tight closure, this gives
$\fra^r\subseteq J^{[q]}$,
a contradiction.

We prove now ${\rm (}6{\rm )}$. If $(\fra+\frb)^r\not\subseteq
J^{[q]}$, then there are $\ell_1$ and $\ell_2$ such that 
$\ell_1+\ell_2=r$ and $\fra^{\ell_1}\not\subseteq J^{[q]}$,
$\frb^{\ell_2}\not\subseteq J^{[q]}$. Therefore 
$\nu_{\fra+\frb}^J(q)\leq\nu^J_{\fra}(q)+\nu_{\frb}^J(q)$
for all $q$, which gives ${\rm (}6{\rm )}$.
\end{proof}

As pointed out in \cite{TW}, the F-pure threshold can be considered
as an analogue of the log canonical threshold. Similarly, 
the F-thresholds play the role of the jumping coefficients from
\cite{ELSV}. We will see more clearly this analogy in the next sections.

In what follows we fix the ideal $\fra$, and study
the F-thresholds which appear for various $J$.
We record in the next proposition some easy properties which deal with
the variation of $J$.

\begin{proposition}\label{prop2}
\begin{enumerate}

\item If $\fra$ and $J_1$, $J_2$ are as above with $J_1\subseteq J_2$,
then $\pt^{J_2}(\fra)\leq\pt^{J_1}(\fra)$.
In particular, the F-pure threshold $\pt(\fra)$ is the smallest
(nonzero) F-threshold of $\fra$.

\item If $J=\bigcap_{\lambda\in\Gamma}J_{\lambda}$, then
$$\pt^J(\fra)=\sup_{\lambda\in\Gamma}\pt^{J_{\lambda}}(\fra).$$

\item We have $\pt^{J^{[q]}}(\fra)=q\cdot\pt^J(\fra)$ for every $q$.
\end{enumerate}
\end{proposition}

\begin{proof}
The first assertion is straightforward, as we have 
$\nu^{J_2}_{\fra}(q)\leq\nu^{J_1}_{\fra}(q)$ for all $q$. For the second 
assertion, note that since the Frobenius morphism is flat, we have
$J^{[q]}=\bigcap_{\lambda}J_{\lambda}^{[q]}$, so 
$\nu^J_{\fra}(q)=\max_{\lambda}\nu^{J_{\lambda}}_{\fra}(q)$
which gives the formula for $\pt^J(\fra)$.
The equality in
${\rm (}3{\rm )}$ is trivial, as in the definition of $\pt^J(\fra)$
we have a limit.
\end{proof}

\bigskip

When the ideal $\fra$ is generated by one element, then we can say
more. The next proposition shows that in this case the F-threshold
determines the numbers $\nu_{\fra}^J(q)$ for all $q$.
If $\fra=(f)$, we simply write $\nu^J_f(q)$ and $\pt^J(f)$.
We denote by $\lceil \alpha\rceil$ the smallest integer $\geq\alpha$.

\begin{proposition}\label{prop3}
Let $J\subseteq\frmm$ be an ideal whose radical contains 
$f\neq 0$. 
For every $q$ we have
$$\frac{\nu^J_f(pq)+1}{pq}\leq\frac{\nu^J_f(q)+1}{q},$$
so $\pt^J(f)=\inf_q\frac{\nu^J_f(q)+1}{q}$.
Moreover, we have $\nu^J_f(q)=\lceil\pt^J(f)q\rceil-1$ for all $q$.
\end{proposition}

\begin{proof}
For the first assertion it is enough to note that
if $f^{\nu_{\fra}^J(q)+1}$ lies in  $J^{[q]}$ then
$f^{p(\nu_f^J(q)+1)}$ is in $J^{[pq]}$.
The last statement follows from 
$$\frac{\nu^{J}_f(q)}{q}<\pt^J(f)\leq\frac{\nu^J_f(q)+1}{q}.$$
\end{proof}

We clearly have $\pt^{(f)}(f)=1$, so $1$ is always an F-threshold
for principal ideals. The next proposition shows that moreover, in this 
case it is enough to understand the thresholds in $(0,1)$.

\begin{proposition}\label{prop4}
If $J$ is an ideal containing the nonzero $f$ in its radical, then
\begin{equation}
\pt^{fJ}(f)=\pt^J(f)+1,\,\,\pt^{(J:f)}(f)=\max\{\pt^J(f)-1,0\}.
\end{equation}
In particular, a nonnegative $\lambda$ is an F-threshold of
$\fra$ if and only if $\lambda+1$ is. 
\end{proposition}

\begin{proof}
The proof is straightforward. The only thing to notice is that
since the Frobenius morphism is flat, we have $(J:f)^{[q]}
=J^{[q]}:f^q$ for all $q$.
\end{proof}

\begin{remark}
It follows easily from Proposition~\ref{prop3}
that when $\fra$ is principal, $\pt^J(\fra)$ is a rational
number if and only if
 the function 
$e\to \overline{\nu}^J_{\fra}(p^e):=
\nu_{\fra}^J(p^{e})-p\nu^J_{\fra}(p^{e-1})$
is eventually periodic. Furthermore, this is equivalent with the
fact that the series
$$P_{\fra}^J(t)=\sum_{e\geq 1}\nu^J_{\fra}(p^e)t^e$$
is a rational function.

One could ask more generally whether for any $\fra$ the
above series is a rational function (again, this would imply
that $\pt^J(\fra)$ is rational).  
It follows from \cite{BMS1}
 that this stronger assertion holds for monomial ideals.
In fact, in this case it is again true that the
function $\overline{\nu}^J_{\fra}$
is eventually periodic.
\end{remark}

\begin{remark}
In the study of singularities in characteristic zero, one can often reduce
the invariant of an arbitrary ideal $\fra$ to that of a principal ideal $(f)$
by taking $f$ general in $\fra$. This does not work in our setting.
For example, let $\fra=\frmm^{[p]}$. We have 
$\pt(\fra)=\pt(\frmm^p)=n/p$, but for every $f\in\fra$,
$\nu_f(p)=0$, so $\pt(f)\leq 1/p$.
\end{remark}

\section{F-thresholds as jumping coefficients}

Test ideals are a very useful tool in tight closure theory.
In \cite{HY} Hara and Yoshida introduced a generalization of test
ideals in the setting of pairs. These ideals enjoy properties similar
to those of multiplier ideals in characteristic zero. In fact, there
is a strong connection between the test ideals and the multiplier ideals
via reduction mod $p$ (see \cite{HY}, and also the next section).
We start by reviewing their definition in our particular setting, in order
to describe the connection between test ideals and F-thresholds.

Let us fix first some notation. Let $E=E(R):=H^n_{\frmm}(R)$ be the top 
local cohomology 
module of $R$, so $E$ is isomorphic to the injective hull of $k$.
If $x_1,\ldots,x_n$ form a regular system of parameters in $R$, then
the completion
$\hat{R}$ of $R$ is isomorphic
to the formal power series ring $k[[X_1,\ldots X_n]]$  
such that $x_i$ corresponds to $X_i$. Note that we have
\begin{equation}\label{isom1}
E(R)\simeq E(\hat{R})\simeq R_{x_1\ldots x_n}/\sum_{i=1}^nR_{x_1\ldots
\widehat{x_i}\ldots x_n}.
\end{equation}
Whenever working in $E$ we will assume we have fixed such a
regular system of parameters, so
via the above isomorphism we may represent each element of $E$ as 
the class $[u/(x_1\ldots x_n)^d]$ for some $u\in R$ and some $d$.  

We will use freely Matlis duality: ${\rm Hom}(-,E)$
induces a duality between finitely generated $\hat{R}$-modules and
Artinian $\hat{R}$-modules (which are the same as the Artinian $R$-modules).
See, for example \cite{BH} for more on local cohomology and Matlis duality.

On $E$ we have a Frobenius morphism $F_E$ which via the isomorphism
in (\ref{isom1}) is given by
$$F_E([u/(x_1\ldots x_n)^d])=[u^p/(x_1\ldots x_n)^{pd}].$$
$F_E$ is injective. Moreover, if $a\in R\setminus \{0\}$ 
is such that $aF_E^e(w)=0$ for all $e$, then $w=0$.
Indeed, this is an immediate consequence of the fact that every
ideal is equal with its tight closure.

Let $\fra\subseteq\frmm$ be a fixed ideal.
For every $r\ge 0$ and $e\geq 1$ we put
$$Z_{r,e}:=\Ker(\fra^rF_E^e)=\{w\in E\vert hF_E^e(w)=0\,{\rm for}\,{\rm all}\,
h\in\fra^r\}.$$

\begin{lemma}\label{inclusion}
If $r<s$, then $Z_{r,e}\subseteq Z_{s,e}$. We have $E=\bigcup_r Z_{r,e}$.
Moreover, $Z_{pr,e+1}$ is contained in $Z_{r,e}$.
\end{lemma}

\begin{proof}
The first assertions are clear, and the last one follows from the
injectivity of $F_E$ and the fact that $F_E(hw)=h^pF_E(w)$ for
all $h\in R$ and $w\in E$. 
\end{proof}

\begin{definition}{\rm (}\cite{HY}{\rm )}
If $\fra\subseteq\frmm$ is a nonzero ideal, and if $c\in\RR_+$,
the test ideal of $\fra$ of exponent $c$ is
$$\tau(\fra^c):={\rm Ann}_R\left(\bigcap_{e\geq 1}Z_{\lceil cp^e\rceil,e}\right).$$
As $E$ is Artinian,
it follows from Lemma~\ref{inclusion} that $\tau(\fra^c)=
{\rm Ann}_RZ_{\lceil cp^e\rceil,e}$
if $e\gg 0$. 
\end{definition}

For every $c>0$, let $Z_c:=\bigcap_eZ_{\lceil cp^e\rceil,e}$.
Note that  $Z_c\neq E$. Indeed, if 
$m\geq c$ is an integer and if $h$ is a nonzero element of $\fra$, 
then $Z_c\subseteq Z_{mp^e,e}\subseteq {\rm ker}(h^{mp^e}F_E^e)$,
which is equal to the kernel of the multiplication by $h^m$ on $E$
(this follows from the
injectivity of $F_E$). Therefore $Z_c$ is a proper submodule of $E$. 

If we replace $R$ by $\hat{R}$ and $\fra$ by $\fra \hat{R}$, then
$Z_c$ remains the same. 
For the basic properties of test ideals we refer the reader to \cite{HY}.
We prove only the following Lemma which we will need in the
next section. See \cite{HY}, \cite{HW} and \cite{Smi} for related 
stronger statements.

\begin{lemma}\label{characterization}
For every $c>0$, the submodule
$Z_c$ is the unique maximal 
proper submodule
of $E$ invariant by all $hF_E^e$, where $e\geq 1$ and
$h\in\fra^{\lceil cp^e\rceil}$.
\end{lemma}

\begin{proof}
It is clear that $Z_c$ is invariant under $hF^e_E$
as above, as $hF_E^e(Z_{\lceil cp^e\rceil,e})=0$
by definition. Since $Z_c$ does not change when we pass
to the completion, we may assume that $R$ is complete.

In this case every proper submodule of $E$ has
nonzero annihilator.
Therefore
in order to finish the proof it is enough to show that 
if $g\in R$ is a nonzero element, and if $w\in E$ is such that
$g\fra^{\lceil cp^e\rceil}F_E^e(w)=0$ for all $e\geq 1$, then
$w\in Z_c$. Fix $e'$ and $h\in\fra^{\lceil cp^{e'}\rceil}$.
For every $e$ we have $gF^e_E(hF_E^{e'}(w))=gh^{p^e}F_E^{e+e'}(w)=0$,
as $h^{p^e}\in\fra^{p^e\lceil cp^{e'}\rceil}\subseteq\fra^{\lceil cp^{e+e'}
\rceil}$. This implies that $hF_E^{e'}(w)=0$, which completes the proof.
\end{proof}

Our goal is to show that the F-thresholds we have introduced in the
previous section can be interpreted as jumping coefficients for
the test ideals. We start by interpreting the function $\nu^J_{\fra}$
in terms of the Frobenius morphism on $E$.

\begin{lemma}\label{interpretation}
Let $\fra$ and $J\subseteq\frmm$ be nonzero ideals, with $\fra$
contained in the radical of $J$. If $M$ is a submodule of $E$
such that $J={\rm Ann}_R(M)$, then
$\nu^J_{\fra}(p^e)$ is the largest $r$ such that $M\not\subseteq Z_{r,e}$.
\end{lemma} 

\begin{proof}
For every $w\in M$ we
put $J_w={\rm Ann}_Rw$, so $J=\bigcap_{w\in M}J_w$. 
If $w=[u/(x_1\ldots x_n)^d]$, then $J_w^{[p^e]}=(x_1^d,\ldots,x_n^d)^{[p^e]}
:u^{p^e}$. For every $w$, we see that
$\nu_{\fra}^{J_w}(p^e)$ is the largest $r$ such that 
$w\not\in Z_{r,e}$. As $\nu^J_{\fra}=\max_{w\in M}\nu^{J_w}_{\fra}$,
we get the assertion in the lemma.
\end{proof}

\begin{remark}
By Matlis duality, we may take in the above Lemma $M={\rm Ann}_E(J)$.
\end{remark}

For future reference we include also the next lemma
whose proof is immediate from definition.

\begin{lemma}\label{inclusion2}
If $\fra\subseteq\frb$, then $\tau(\fra^c)\subseteq\tau(\frb^c)$
for every $c\in\RR_+$. If $c_1<c_2$, then $\tau(\fra^{c_2})\subseteq
\tau(\fra^{c_1})$ for every ideal $\fra$.
\end{lemma}

\smallskip

\begin{proposition}\label{inclusion3}
If $\fra\subseteq\frmm$ is a nonzero ideal contained in the radical of $J$,
then
$$\tau(\fra^{\pt^J(\fra)})\subseteq J.$$
Going the other way, if $\alpha\in\RR_+$,
then $\fra$ is contained in the radical of $\tau(\fra^{\alpha})$ and
$$\pt^{\tau(\fra^{\alpha})}(\fra)\leq \alpha.$$
Therefore the maps $J\longrightarrow\pt^J(\fra)$ and $\alpha\longrightarrow
\tau(\fra^{\alpha})$
give a bijection between the set of test ideals of $\fra$
and the set of F-thresholds of $\fra$.
\end{proposition}

\begin{proof}
For the first statement, let $M={\rm Ann}_EJ$, so by Matlis duality we need
to prove that $M\subseteq Z_{\lceil\pt^J(\fra)p^e\rceil,e}$ for all $e$.
This follows from Lemma~\ref{interpretation} and the fact that
$\lceil c^J(\fra)p^e\rceil>\nu^J_{\fra}(p^e)$.

We show now that for $\alpha\in\RR_+$, we have 
$\fra\subseteq {\rm Rad}(\tau(\fra^{\alpha}))$. 
Let $e\gg 0$ be such that $\tau(\fra^{\alpha})
={\rm Ann}_RZ_{\lceil\alpha p^e\rceil,e}$. 
If $m\geq\alpha$ is an integer, it follows from the
injectivity of $F_E$ that $\fra^m\subseteq\tau(\fra^{\alpha})$.

We deduce now from Lemma~\ref{interpretation}
and from the definition of $\tau(\fra^{\alpha})$
that 
$$\nu_{\fra}^{\tau(\fra^{\alpha})}(p^e)\leq\lceil\alpha p^e\rceil -1
<\alpha p^e.$$
Dividing by $p^e$ and taking the limit gives the required inequality.

The last statement is a formal consequence of the first two assertions.
\end{proof}

\begin{remark}\label{minimal}
It follows from Proposition~\ref{inclusion3}
that if we have an F-threshold $c$ of $\fra$,
then there is a unique minimal ideal $J$ such that $\pt^J(\fra)=c$.
Indeed, this is $\tau(\fra^c)$. Moreover, if $c_1$ and $c_2$
are such F-thresholds, then $c_1<c_2$ if and only if $\tau(\fra^{c_2})$
is strictly contained in $\tau(\fra^{c_1})$.  
\end{remark}

\begin{remark}
As $R$ is Noetherian, it follows from the previous remark that
there is no strictly decreasing sequence of F-thresholds of $\fra$. 
\end{remark}

\begin{remark}
There are arbitrarily large thresholds: take, for example
$p^e\pt(\fra)=\pt^{\frmm^{[p^e]}}(\fra)$ 
for $e\geq 1$. Note also that a sequence of thresholds
$\{c_m\}_m$ of $\fra$ 
is unbounded if and only if $\bigcap_m\tau(\fra^{c_m})=(0)$.
The only thing 
one needs to check is that $\bigcap_{c\in\RR_+}\tau(\fra^c)=(0)$.
This follows since for every integer $\ell\geq n$, we have
$\tau(\fra^{\ell})\subseteq\tau(\frmm^{\ell})\subseteq\frmm^{\ell-n+1}$.
\end{remark}

The jumping coefficients for multiplier ideals are discrete. We do not
know if the analogous assertion is true for the F-thresholds.

\begin{question}\label{accumulation}
Given an ideal 
$(0)\neq\fra\subseteq\frmm$, could there exist finite accumulation
points for the set of F-thresholds of $\fra$ ?
\end{question}

\begin{remark}
Given a test ideal $J$ corresponding to $\fra$, the set of those
$\alpha\in\RR_+$ such that $\tau(\fra^{\alpha})=J$ is an interval
of the form $[a,b)$. Indeed, if $a=\pt^J(\fra)$, it follows from
Lemma~\ref{inclusion2} and Proposition~\ref{inclusion3}
that for $\lambda<a$, 
$\tau(\fra^{\lambda})$ strictly contains $J=\tau(\fra^a)$.
On the other hand, if
$J'$ is the largest 
test ideal strictly contained in $J$, and if $b=\pt^{J'}(\fra)$,
then it is clear that $b=\sup\{\lambda\vert\tau(\fra^{\lambda})=J\}$,
which gives our assertion. 
\end{remark}

\begin{example}\label{monomial}
Consider the case when $\fra$ is a monomial ideal, i.e.
$\fra$ is  generated by monomials in the localization of $k[X_1,\ldots,X_n]$
at $(X_1,\ldots,X_n)$. It is shown in \cite{HY}
that for every $\alpha$ we have $\tau(\fra^{\alpha})
={\mathcal I}(\fra^{\alpha})$, where ${\mathcal I}(\fra^{\alpha})$ 
is the multiplier ideal of $\fra$ with exponent
$\alpha$. 
It follows from this and from Proposition~\ref{inclusion3}
that the set of F-thresholds of $\fra$ coincides with the
set of jumping coefficients of the multiplier ideals of $\fra$.

Let us recall the description of multiplier ideals for monomial ideals
from \cite{howald}. Consider the Newton polyhedron $P_{\fra}$ of $\fra$:
this is the convex hull in ${\mathbb R}^n$ of 
$\{u\in\NN^n\vert X^u\in\fra\}$,
where for $u=(u_1,\ldots,u_n)$ we put $X^u=X_1^{u_1}\ldots X_n^{u_n}$. 
If we put $e=(1,\ldots,1)$, then 
\begin{equation}
{\mathcal I}(\fra^{\alpha})=
(X^u\vert u+e\in {\rm Int}(\alpha\cdot P_{\fra})).
\end{equation}
It follows that each jumping coefficients $\alpha$ of the multiplier ideals
is associated to some $b=(b_i)$ with all $b_i$ positive integers,
where $\alpha$ 
is such that $b$ lies in the boundary of $\alpha\cdot P_{\fra}$.
Of course, several distinct $b$ can give the same $\alpha$.

In fact, one can show that in order to compute all the F-thresholds
of $\fra$ it is enough to consider only ideals $J$ of the form
$(X_1^{b_1},\ldots,X_n^{b_n})$ with $b_i$ positive integers.
Moreover, one can check directly that if $J$ is the above ideal,
then $\pt^J(\fra)=\alpha$, where $\alpha$ is 
associated to 
$b=(b_1,\ldots,b_n)$ as above (see \cite{BMS1} for this approach).  
\end{example}

\medskip

We end this section by considering in more detail the case
when $\fra=(f)$ is a principal ideal. One can easily check that
for such $\fra$ we have 
$Z_{r,e}=Z_{pr,e+1}$. This shows that any two $Z_{r_1,e_1}$
and $Z_{r_2,e_2}$ are comparable. As $E$ is Artinian, we may consider
the submodules of $E$ defined inductively as follows: let $M_0:=\{0\}$
be the minimal module in ${\mathcal M}:=\{Z_{r,e}\vert r,e\}$, 
and for $m\geq 1$,
let $M_m$ be the unique minimal module in 
${\mathcal M}\smallsetminus\{M_0,\ldots,M_{m-1}\}$.
It follows that $M_{m}$ is properly contained in $M_{m+1}$.
In addition, given any $r$ and $e$, 
 either $M_m\subseteq Z_{r,e}$,
or $Z_{r,e}\in\{M_0,\ldots,M_{m-1}\}$.

\begin{proposition}\label{notation}
With $f\in\frmm$ as above, we put
for every $i$,  $J_i={\rm Ann}_R(M_i)$, and  $c_i=\pt^{J_i}(f)$.
\begin{enumerate}
\item For $i\geq 1$, let
 $\nu_i(e)$ be the largest $r$ such that $Z_{r,e}\subseteq M_{i-1}$.
Then 
$$c_i=\lim_{e\to\infty}\frac{\nu_i(e)}{p^e}.$$
\item Every $J_i$ is a test ideal of $\fra$, and if 
$J$ is any test ideal different from all $J_i$, then $J$ is contained
in all these ideals.
\end{enumerate}
\end{proposition}

\begin{proof}
Note that by definition 
$\nu_i(e)$ is the largest $r$ such that $M_i\not\subseteq
Z_{r,e}$. By Lemma~\ref{interpretation}, we get $\nu_i(e)=\nu^{J_i}_{f}(p^e)$
and this proves ${\rm (}1{\rm )}$.

We show that $J_i$ is a test ideal by proving that $\tau(f^{c_i})=J_i$.
By Lemma~\ref{inclusion3}, it is enough to show that 
$J_i\subseteq\tau(f^{c_i})$ for $i\geq 1$. 
This follows 
from $\lceil c_ip^e\rceil=\nu_i(p^e)+1$ (see Proposition~\ref{prop3})
which implies $\bigcap_eZ_{\lceil c_ip^e\rceil,e}=M_i$.

For the last statement it is enough to show that
for all $i\geq 0$ and for $c\in [c_i,c_{i+1})$ we have $\tau(f^c)=J_i$
(with the convention $c_0=0$).
If $e\gg 0$, we have $\lceil cp^e\rceil<\nu_{i+1}(e)$, so
$\bigcap_eZ_{\lceil cp^e\rceil,e}\subseteq M_i$. This implies
$J_i\subseteq \tau(f^c)$, and the other inclusion is clear, 
as we have seen that $\tau(f^{c_i})=J_i$.
\end{proof}

\begin{remark}
Note that in the case of a principal ideal 
the set of F-thresholds of $(f)$ is discrete 
if and only if $\lim_{m\to\infty}c_m=\infty$.
This is equivalent with the fact that $\bigcup_mZ_m=E$.
Note also that by 
the periodicity of the F-thresholds (see Lemma~\ref{prop4}),
this is further equivalent with the finiteness of the 
set of F-thresholds in $(0,1)$.
\end{remark}

It follows from Proposition~\ref{prop3}
that $\pt(f)=1$ if and only if 
$\nu_f(p^e)=p^e-1$ for every $e$. The following proposition
based on an argument of Fedder
shows that in fact, it is enough to check this for only one $e\geq 1$.

\begin{proposition}{\rm (}\cite{Fe}{\rm )}\label{one_value}
If $f$ is a nonzero element in $\frmm$, then
$\pt(f)=1$ if and only if there is $e$ such that $\nu_f(p^e)=p^e-1$.
Moreover, this is the case if and only if the action of the Frobenius
morphism on $H^{n-1}_{\frmm}(R/(f))$ is injective.
\end{proposition}

\begin{proof}
The exact sequence 
$$0\to R\to R\to R/(f)\to 0$$
induces an isomorphism of $H^{n-1}_{\frmm}(R/(f))$ with
the annihilator of $f$ in $E$. Moreover, 
via this identification the Frobenius morphism
on $H^{n-1}_{\frmm}(R/(f))$ is given by $\widetilde{F}(u)=f^{p-1}F_E(u)$.

We see that $\widetilde{F}^e$ is injective if and only if 
$Z_{p^e-1,p^e}=(0)$. This is the case if and only if $\nu_f(p^e)=p^e-1$.
Since $\widetilde{F}$ is injective if and only if $\widetilde{F}^e$ is,
this completes the proof.
\end{proof}

\section{Reduction mod $p$ and the connection with the Bernstein polynomial}

In this section
 we study the way our invariants behave for different reductions
mod $p$ of a given ideal.
Everything in this section works in the usual framework
for reducing mod $p$ which is used in tight closure theory
(see for example \cite{HY}). In order to simplify the presentation
as well as the notation, we prefer to work in the following 
concrete setup. The interested reader should have no trouble 
translating everything to the general setting.

Let $A$ be the localization of $\ZZ$ at some nonzero integer.
We fix a nonzero ideal $\fra$ 
of $A[X]=A[X_1,\ldots,X_n]$, such that
$\fra\subseteq (X_1,\ldots,X_n)$. Let ${\mathbb F}_p=\ZZ/p\ZZ$. 
We want to relate the invariants 
attached to $\fra_{\QQ}:=\fra\cdot\QQ[X]$ 
around the origin with those attached to
the localizations of the reductions mod $p$,
$$\fra_p:=\fra\cdot {\mathbb F}_p[X]_{(X_1,\ldots,X_n)},$$
where $p$ is a large prime.
We will use the same subscripts 
whenever tensoring with $\QQ$ or reducing (and localizing) mod $p$.
Note that
since we are interested only in large primes, we are free
to further localize $A$ at any nonzero element.

Let us consider a log resolution of
$\fra_{\QQ}$ 
defined over $\QQ$: this is a proper birational morphism
$\pi_{\QQ} : Y_{\QQ}\longrightarrow {\mathbb A}^n_{\QQ}$, 
with $Y_{\QQ}$ 
smooth,
 such that the product between
$\pi^{-1}_{\QQ}(\fra_{\QQ})$ and the
ideal defining the exceptional locus of $\pi_{\QQ}$ is principal, and it defines
a divisor with simple normal crossings. Such a resolution exists 
by \cite{hironaka}.
After further localizing $A$, we may assume that $\pi_{\QQ}$ is obtained
by extending the scalars from a morphism $\pi : Y\longrightarrow 
{\mathbb \A}_A^n$
with analogous properties.

If we denote by $D$ the effective divisor defined by 
 $\pi^{-1}(\fra)$ and if $K$ 
is the relative canonical divisor of $\pi$
(i.e. the effective divisor defined by the Jacobian of $\pi$), then
for all $\alpha\in\RR$ 
\begin{equation}\label{mult1}
\cI(\fra^{\alpha}):=H^0(Y,\cO_{Y}(K-\lfloor\alpha D\rfloor)).
\end{equation}
Here $\lfloor\alpha D\rfloor$
denotes the integral part of the $\RR$-divisor $\alpha D$.
Note that $\cI(\fra^{\alpha})_{\QQ}$ is the multiplier ideal of
$\fra_{\QQ}$ of with exponent $\alpha$.
We refer for the theory of multiplier ideals to \cite{lazarsfeld}.

The jumping 
numbers at $(0,\ldots,0)$ introduced in \cite{ELSV} are the numbers $\lambda$
such that $\cI(\fra^{\lambda})_{\QQ}$ is strictly contained in 
$\cI(\fra^{\lambda-\epsilon})_{\QQ}$ in a neighborhood of the origin,
 for every $\epsilon>0$.
The smallest positive such number is the log canonical threshold
$\lct_0(\fra)$: it is the first $\lambda$ such that 
$\cI(\fra^{\lambda})_{\QQ}$ is different from the structure sheaf
around the origin. In order to simplify the notation we will drop 
the subscript $\QQ$ whenever considering the invariants associated to $\fra_{\QQ}$.

In our setting, by taking $p\gg 0$, we may assume that the above 
resolution induces a log resolution $\pi_p$ for $\fra_p$. Over $\QQ$
we have $R^i\pi_*(\cO_Y)_{\QQ}=0$ for $i\geq 1$. This remains true for
the reductions mod $p$ if $p\gg 0$. From now on we assume that
$p$ is large enough, so these conditions are satisfied.
We define $\cI(\fra_p^{\alpha})$ by a formula similar to (\ref{mult1}),
using $\pi_p$. Note that for a fixed $\alpha$, we have 
$\cI(\fra^{\alpha})_p=\cI(\fra_p^{\alpha})$
if $p\gg 0$.

We recall two results which 
describe what is known about the connection between multiplier ideals
and test ideals. The first one is proved in more generality in
\cite{HY},
based on ideas from \cite{HW}. 
We include the proof as it is quite short in our context.

\begin{theorem}\label{thm1}
With the above notation, if $p\gg 0$, then for every $\alpha$ we have
$$\tau(\fra_p^{\alpha})\subseteq\cI(\fra_p^{\alpha}).$$
\end{theorem}

\begin{proof}
Let $R$ be the localization of
${\mathbb F}_p[X]$ at $(X_1,\ldots,X_n)$, and let $\frmm$ be
its maximal ideal. 
We denote by $W\subseteq Y_p$ the subset defined by
$\pi_p^{-1}(\frmm)$.
We will use the notation from the previous section.

If $E=H^n_{\frmm}(R)$, using the fact that
 the higher direct images of ${\mathcal O}_{Y_p}$
are zero, 
and the long exact sequence for local cohomology we get
$E\simeq H^n_W(\cO_{Y_p})$. 
A version of Local Duality shows that if
$$\delta : H^n_W(\cO_{Y_p})\longrightarrow H^n_W(\cO_{Y_p}(\lfloor
\alpha D_p\rfloor))$$
is the surjective morphism induced by the natural inclusion of sheaves, then
$\cI(\fra_p^{\alpha})={\rm Ann}_R(\ker\,\delta)$.
By Lemma~\ref{characterization}, it is
therefore enough to show that if $h\in\fra_p^{\lceil \alpha p^e\rceil}$,
then $hF_E^e(\ker\,\delta)\subseteq\ker\,\delta$. 

The Frobenius morphism on local cohomology is induced
by the Frobenius morphism $F$ on the fraction field of $R$.
As the inclusion $hF^e(\cO_{Y_p})\subseteq\cO_{Y_p}$ is clear, in order
to finish it is enough to show also that 
$hF^e(\cO_{Y_p}(\lfloor\alpha D_p\rfloor))
\subseteq\cO_{Y_p}(\lfloor\alpha D_p\rfloor)$.
This is an immediate consequence of the definitions.
\end{proof}

The proof of the next Theorem is more involved, so we refer the reader
to \cite{HY}.

\begin{theorem}\label{thm2}
With the above notation, if $\alpha$ is given
and if $p\gg 0$ (depending on $\alpha$), then
$$\tau(\fra_p^{\alpha})=\cI(\fra^{\alpha}_p).$$
\end{theorem}

We reformulate the above results in terms of thresholds.
In order to do this we index the jumping coefficients of $\fra_{\QQ}$
at the origin by analogy with the F-thresholds, as follows. 
Suppose that $J\subseteq (X_1,\ldots, X_n)A[X]$ is an ideal
containing $\fra$ in its radical. We define
$$\lambda^J_0(\fra):=
\min\{\alpha>0\,\vert\,\cI(\fra^{\alpha})_{\QQ}\subseteq J_{\QQ}
\,{\rm around}\,0\}.$$
It is clear that this is a jumping coefficient of $\fra_{\QQ}$ around the
origin, and 
that every such coefficient appears in this way for a suitable $J$.
For example,  if $J=(X_1,\ldots,X_n)$, then
$\lambda^{J}_0(\fra)=\lct_0(\fra)$.

Using Proposition~\ref{inclusion3}, we may reformulate the above results
as follows. We will denote the invariants of $\fra_p$ with respect to $J_p$,
which we have introduced in \S 1, 
simply by $\pt^J(\fra_p)$ and $\nu^J_{\fra}(p^e)$.

\begin{theorem}\label{thm1'}
If $p\gg 0$, then for every ideal $J$ as above we have
$\pt^J(\fra_p)\leq\lambda^J_0(\fra)$.
In particular, we have the following inequality between the
F-pure threshold and the log canonical threshold:
$\pt(\fra_p)\leq\lct_0(\fra)$.
\end{theorem}

\begin{theorem}\label{thm2'}
Given an ideal $J$ as above, we have
$$\lim_{p\to\infty}\pt^J(\fra_p)=\lambda^J_0(\fra).$$
In particular, we have $\lim_{p\to\infty}\pt(\fra_p)=\lct_0(\fra)$.
\end{theorem}

\begin{remark}
The fact that in Theorem~\ref{thm2} $p$ depends on $\alpha$ is reflected
in Theorem~\ref{thm2'} in that we may have 
$\pt^J(\fra_p)<\lambda^J_0(\fra)$ for infinitely many $p$.
This is a very important point, 
and we will see 
examples of such a behavior (for $J=\frmm$) in the next section.
\end{remark}

We discuss now possible further connections 
between the
invariants over $\QQ$ and those of the reductions mod $p$.
We formulate 
them in the case $J=\frmm$
and we will give some examples in \S 4.
However, note that similar questions can be asked for arbitrary $J$.

\begin{conjecture}\label{infinite}
Given the ideal $\fra$, there are infinitely many primes $p$
such that $\pt(\fra_p)=\lct_0(\fra)$.
\end{conjecture}

\begin{problem}\label{cong1}
Given the ideal $\fra$, give conditions such that
there is a positive integer $N$ with the following property:
for every prime $p$ with
$p\equiv 1$ (mod $N$) we have $\pt(\fra_p)
=\lct_0(\fra)$. 
\end{problem}

\begin{problem}\label{limit_p}
Give conditions on an ideal $\fra$
such that there is a positive integer $N$, and
rational functions $R_i\in\QQ(t)$
for every $i\in\{0,\ldots,N-1\}$ with ${\rm gcd}(i,N)=1$
with the following property:
$\pt(\fra_p)=R_i(p)$ whenever $p\equiv i$ (mod $N$) and
$p$ is large enough.
\end{problem}

These problems are motivated by the examples we will discuss
in the next section. We will see that the behavior described
in the problems is satisfied in many cases. On the other hand,
Example~\ref{elliptic} below shows that one can not expect for
such a behavior to hold in general. We will see in this example
that the failure is related to subtle arithmetic phenomena.

However, note that if $p$ is an odd prime, then one can reinterpret 
the condition $p\equiv 1$ (mod $N$) 
in Problem~\ref{cong1} as saying that $p$ is completely split
in the cyclotomic field of the $N$th roots of unity
(see \cite{Neukirch}, Cor.~10.4). We will see
that somethings similar happens in Example~\ref{elliptic} below: 
there is a number field $K$ such that if $p$ splits completely in $K$,
then the log canonical threshold is equal to the corresponding 
F-pure threshold. This motivates the following

\begin{question}\label{number_field}
Given an ideal $\fra$ as above, is there a number field $K$
such that whenever the prime $p\gg 0$ splits completely in $K$, we have 
$\pt(\fra_p)=\lct(\fra)$ ?
\end{question}

Note that by \v{C}ebotarev's Density Theorem
(see \cite{Neukirch}, Cor. 13.6), given a number field $K$
there are infinitely many primes $p$ which split completely in $K$.
Therefore a positive answer to Question~\ref{number_field}
would imply Conjecture~\ref{infinite}.

\medskip

We include here another problem with a similar flavor, on the behavior of 
the functions $\nu^{J}_{\fra}(p^e)$ when we vary $p$. The interest in
this problem comes from the fact that whenever we can prove that 
such a behavior holds,
one can use this 
to give roots of the Bernstein-Sato polynomial of $\fra_{\QQ}$
(see Remark~\ref{roots} below). The Conjecture is proved for monomial
ideals in \cite{BMS1}. For other examples, see the next section.

\begin{problem}\label{classes}
Find conditions on an ideal $\fra$ such that the following holds.
Given an ideal $J$ as above, and $e\geq 1$,
there is a positive integer $N$, 
and polynomials $P_j\in\QQ[t]$ of degree $e$,
for every $j\in\{1,\ldots,N-1\}$ with ${\rm gcd}(j,N)=1$,
such that $\nu^{J}_{\fra}(p^e)=P_{j}(p)$
for every $p\gg 0$, $p\equiv j$ (mod $N$). 
When could $N$ be chosen independently on
$J$ and $e$ ?
\end{problem}

\bigskip

We turn now to a different connection between invariants
which appear in characteristic zero and the ones
we have defined in \S 1. The characteristic zero
invariants we will consider
are the roots of 
the Bernstein-Sato polynomial, whose definition we now recall.

Let $I\subseteq\CC[X_1,\ldots,X_n]$ be a nonzero ideal, and
let $f_1,\ldots,f_r$ be nonzero generators of $I$. 
We introduce indeterminates
$s_1,\ldots,s_r$ and the Bernstein-Sato polynomial
$b_I$ is the monic polynomial in one variable of minimal degree
such that we have an equation
\begin{equation}\label{eq20}
b_{I}(s_1+\ldots+s_r)\prod_{i=1}^rf_i^{s_i}
=\sum_cP_c(s,X,\partial_X)
\bullet\prod_{j,c_j<0}{s_j\choose -c_j}\prod_{i=1}^rf_i^{s_i+c_i}.
\end{equation}
Here the sum varies over finitely many $c\in\ZZ^r$ such that
$\sum_jc_j=1$, for every such $c$ we have the nonzero
differential operator
$P_c\in\CC[s_j, X_i,\partial_{X_i}\vert j\leq r,i\leq n]$,
and as usual ${s_j\choose {-c_j}}=s_j(s_j-1)\ldots(s_j+c_j+1)/(-c_j)!$.
Note that $\bullet$ denotes the action of a differential operator.
Equation (\ref{eq20}) is understood formally, but if we let
$s_i=m_i\in\NN$, then it has the obvious meaning.
If we require $(\ref{eq20})$ to hold only in some neighborhood of 
the origin in $\CC^n$, 
then we get the local Bernstein-Sato polynomial $b_{I,0}(s)$. 

Note that if $r=1$, i.e. if $I=(f)$ is a principal
ideal, then (\ref{eq20}) takes the more familiar form
$$b_f(s)f^s=P(s,X,\partial_X)\bullet f^{s+1}.$$
We refer to \cite{Bj} for some basic properties properties
of the Bernstein-Sato polynomial of principal ideals, and
to \cite{BMS2} for the general case.

In the case of principal ideals,
there is an extensive literature on connections between this polynomial 
and other invariants of singularities (see
\cite{malgrange}, \cite{kashiwara}, \cite{igusa} and \cite{kollar}).
Some of these results have been extended to arbitrary ideals in
\cite{BMS2}.

Here are a few properties which are relevant to our study.
First, it is proved in \cite{BMS2} that this polynomial does not
depend on the choice of generators. All the roots 
of $b_{I,0}$ are negative
rational numbers, the largest one is $-\lct_0(I)$, and 
for every jumping coefficient around the origin $\lambda$ of $I$,
if $\lambda\in [\lct_0(I),\lct_0(I)+1)$, then $-\lambda$
is a root 
of $b_{I,0}$.
For these facts, see
 \cite{kashiwara},
\cite{kollar} and \cite{ELSV} for principal ideals and
\cite{BMS2} for the general case.

\smallskip

We return now to our setting.
The extension of our ideal $\fra$ to $\CC$ defines
a Bernstein-Sato polynomial around the origin, which we simply denote
by $b_{\fra,0}$. Consider the defining equation (\ref{eq20})
and let $B$ be a subalgebra of $\CC$, finitely generated over $\ZZ$
and containing all the coefficients of $b_{\fra,0}$ and of the $P_c$.
Moreover, we may assume that for all $c$ which appear in (\ref{eq20})
and for all $j$ such that $c_j<0$, $(-c_j)!$ is invertible in $B$.

It is clear that there is $M$ such that for every prime $p\geq M$,
there is a maximal ideal $P$ of $B$ with $pA=P\cap A$.
For such $p$ and $P$, let $R_p$ and $S_P$ be the localizations
of ${\mathbb F}_p[X]$ and 
$(B/P)[X]$, respectively, at the ideal generated by the variables.
Suppose now that $J\subseteq 
(X_1,\ldots,X_n)A[X]$ is an ideal containing $\fra$
in its radical. We will denote by $J_p$ and $J_P$ the image of
$J$ in $R_{p}$ and $S_P$, respectively. 

Note that since $S_P$ is flat over $R_p$ and since the Frobenius morphism 
is flat,
it follows that for every $e$ we have 
$J_P^{[p^e]}\cap R_{p}=J_p^{[p^e]}$. In particular, we have
\begin{equation}\label{equality_for_nu}
\nu^{J}_{\fra}(p^e)=\nu^{J_P}_{\fra_pS_P}(p^e).
\end{equation}

\begin{proposition}\label{bernstein}
If $\fra\subseteq (X_1,\ldots,X_n)A[X]$ is a nonzero ideal,
then for every prime $p\gg 0$ and for every $J\subseteq(X_1,\ldots,X_n)A[X]$
containing $\fra$ in its radical we have
\begin{equation}\label{eq22}
b_{\fra,0}(\nu^{J}_{\fra}(p^e))=0\,{\rm in}\,{\mathbb F}_p
\end{equation}
for all $e$.
\end{proposition}

\begin{remark}
Recall that all roots of $b_{\fra,0}$ are rational, so 
$b_{\fra,0}\in\QQ[s]$. After localizing $A$ at a suitable element,
we may assume that $b_{\fra,0}\in A[s]$.
Therefore for every $m\in\ZZ$,  $b_{\fra,0}(m)$ has a well-defined class
in ${\mathbb F}_p$.
\end{remark}

\begin{proof}[Proof of Proposition~\ref{bernstein}]
We use the above notation and
let $m=\nu^{J}_{\fra}(p^e)$. 
Recall that we have  
generators $f_1,\ldots,f_r$ of $\fra$.
It follows from (\ref{equality_for_nu})
that there are nonnegative integers
$\ell_1,\ldots,\ell_r$ such that $\sum_i\ell_i=m$
and $\prod_if_i^{\ell_i}\not\in J_P^{[p^e]}$. 
On the other hand, for every nonnegative
integers $\ell'_1,\ldots,\ell'_r$
with $\sum_i\ell'_i=m+1$, we have $\prod_if_i^{\ell'_i}\in J_P^{[p^e]}$.

Note that
(\ref{eq20}) holds in $S_P$ if $s_i=\ell_i$ for all $i$. If
 $c$ and $i$ are such that $\ell_i+c_i<0$, then $\ell_i(\ell_i-1)\ldots
(\ell_i+c_i+1)=0$,
so this term does not appear in the corresponding equality.
As $J_P^{[p^e]}$ is invariant under the action of operators
in $B/P[X_i,\partial_{X_i}\mid i\leq n]$, 
we deduce that $b_{\fra,0}(m)$ is zero in $R_p$, hence in ${\mathbb F}_p$.
\end{proof}

\begin{remark}\label{roots}
Note that whenever we can show that $\fra$ behaves as in
Problem~\ref{classes}, we get roots of $b_{\fra,0}$.
More precisely, suppose that for some $J$ as above and for some $e$,
there is a positive number $N$ and polynomials $P_j\in\QQ[t]$
for every $j$ with ${\rm gcd}(j,N)=1$ such that
$\nu^{J}_{\fra}(p^e)=P_j(p)$
for every prime $p\gg 0$ with $p\equiv j$ (mod $N$).  
In this case, the above proposition 
shows that $b_{\fra,0}(P_j(p))$ is divisible by $p$,
so $p$ divides $b_{\fra,0}(P_j(0))$.
By the Dirichlet Theorem
there are infinitely many such $p$, and therefore $P_j(0)$
is a root of $b_{\fra,0}$.
\end{remark}

\begin{remark}\label{roots_from_coeff}
Let $\fra$ be a principal ideal generated by $f$.
Suppose that the analogue of 
the setup in Problem~\ref{cong1} holds
for a jumping coefficient $\mu\in (0,1]$ (around the origin)
of $f_{\QQ}$. More precisely, suppose that $J$ and $N$ are such that
if $p\equiv 1$ (mod $N$), then $\pt^J(f_p)=\mu$
(a natural choice for such a $J$ is $J={\mathcal I}(f^{\mu})$).
We may choose such $p$ so that $\mu(p-1)$ is an integer.
Since 
$\mu\leq 1$, it follows from Proposition~\ref{prop3}
that in this case we have $\nu^J_{\fra}(p^e)=\mu(p^e-1)$ for all $e\geq 1$.
Remark~\ref{roots} implies now that $-\mu$ is a root of $b_{f,0}(s)$.
As we have already mentioned, this is proved in \cite{ELSV}, but
this would provide an ``explanation'' from our point of view.
\end{remark}

\begin{remark}\label{all}
It is an interesting question which
roots of $b_{\fra,0}$ can be given by the procedure
in Remark~\ref{roots}.
It is proved in \cite{BMS1} 
that this is the case for all the roots 
if $\fra$ is a monomial ideal. On the other case,  
Example~\ref{quadric} below shows that some roots may not come 
from our approach.
\end{remark}

\section{Examples}

\begin{example}\label{quadric}
Let $n\geq 3$ and $f=X_1X_2+X_3^2+\ldots+X_n^2$, so its Bernstein-Sato
polynomial is given by $b_f(s)=(s+1)(s+\frac{n}{2})$
(see \cite{kashiwara1}, Example~6.19, but this is actually one
of the few examples which can be computed directly). We will see
that we can not account for the root $-\frac{n}{2}$
by the procedure described in Remark~\ref{roots}.

We claim that for every $p$ and for every $e\geq 1$,
we have $\nu_f(p^e)=p^e-1$.
To see this note that if over ${\mathbb F}_p$
we have $f^r\in (X_1^{p^e},\ldots,X_n^{p^e})$, then
$(X_1X_2)^r\in (X_1^{p^e},\ldots,X_n^{p^e})$,
as follows by choosing a monomial order on the polynomial ring such that
${\rm in}(f)=X_1X_2$
(see \cite{Eis}, Chapter~15 for monomial orders). 
Therefore $\nu_{f}(p^e)\geq p^e-1$,
so we must have equality.

This shows that the smallest nonzero
F-threshold is $\pt(f_p)=1$. Proposition~\ref{prop4}
shows that if $\lambda$ is an F-threshold of $f_p$ which is
not an integer, then the fractional part of $\lambda$ gives
an F-threshold in $(0,1)$, a contradiction. Therefore
the set of F-thresholds of $f_p$ 
consists of the set of positive integers.

Using Proposition~\ref{prop3} it follows that for every ideal $J$
contained in $(X_1,\ldots,X_n)$ and such that $f\in {\rm Rad}(J)$,
and for every $p$, there is a positive integer $m$ such that
$\nu^J_{f_p}(p^e)=mp^e-1$ for all $e$. Therefore the only root of 
$b_f(s)$ we get by the procedure described in Remark~\ref{roots}
is $-1$.
\end{example}

\begin{example}\label{nondegenerate}
Consider $f\in\ZZ[X_1,\ldots,X_n]$ which we write 
as $f=\sum_{i=1}^rc_iX^{\alpha_i}$,
where all $\alpha_i=(\alpha_{i,1},\ldots,\alpha_{i,n})\in\NN^n$
and all $c_i\in k$ are nonzero. We assume that $\alpha_1,\ldots,\alpha_r$
are affinely independent, 
i.e. if $\sum_i\lambda_i\alpha_i=0$ and $\sum_i\lambda_i
=0$ for $\lambda=(\lambda_i)\in\QQ^r$, then $\lambda_i=0$ for all $i$. 
We will assume also that for every $j\leq n$ there is $i\leq r$ with 
$\alpha_{i,j}>0$ (otherwise we may work in a smaller polynomial ring).

Let $\fra=(X^{\alpha_i}\vert 1\leq i\leq r)$.
One can check that our condition of $f$ implies that
$f$ is generic with respect to $\fra$ in the following sense.
If $P'$ is a compact face of the convex hull $P_{\fra}$ of $\{\alpha_i\vert
1\leq i\leq r\}$, and if $g$ 
is the sum of the terms in $f$ which correspond to
elements in $P'$, then the differential $dg$ does not vanish on $(\CC^*)^n$.
Under this assumption it is proved in \cite{howald1} that around the
origin we have 
${\mathcal I}(\fra^{\alpha})={\mathcal I}(f^{\alpha})$ for $\alpha<1$. 
In particular, we have $\lct_0(f)=\min\{1,\lct(\fra)\}$
(note that $\lct(\fra)=\lct_0(\fra)$ as $\fra$ is a monomial ideal).  

We start by computing $\nu_f(p)$ when $p\gg 0$. 
Over $\ZZ$ we have 
\begin{equation}\label{over_Z}
f^m=\sum_{a}\frac{m!}{a_1!\ldots a_r!}
c_1^{a_1}\ldots c_r^{a_r} X^{\sum_ia_i\alpha_i},
\end{equation}
where the sum is over those $a=(a_i)\in\NN^r$ with $\sum_ia_i=m$.
Our hypothesis on $f$ implies that if $a\neq b$ and $\sum_ia_i=\sum_ib_i$,
then we have $X^{\sum_ia_i\alpha_i}\neq X^{\sum_ib_i\alpha_i}$.

We may assume that $p$ does not divide any of the $c_i$,
so if $m\leq p-1$, then $p$ does not divide any of the coefficients
in (\ref{over_Z}). 
Hence a monomial $X^b$
appears in $f_p^m$ if and only if there are
$a_1,\ldots,a_r\in\NN$ such that $\sum_ia_i=m$ and $\sum_ia_i\alpha_i=b$.
Therefore $f_p^m\in (X_1^p,\ldots,X_n^p)$ if and only if 
the following holds: for every
 $(a_i)\in\NN^r$ with $\sum_ia_i\alpha_{ij}\leq p-1$ for all $j$, we have
$\sum_ia_i\leq m-1$.

Let
$$Q:=\{(a_1,\ldots,a_r)\in\RR_+^r\vert \sum_{i=1}^r a_i\alpha_{ij}\leq 1
\,{\rm for}\,{\rm all}\,j\}.$$
The above discussion shows that 
\begin{equation}\label{eq30}
\nu_f(p)=\min\{p-1,\max_{b\in (p-1)Q\cap\NN^r}\sum_ib_i\}.
\end{equation}

Compare this with the following formula for $\lct(\fra)$
which follows easily from the description of $\lct(\fra)$ in \cite{howald}
(see, for example, Proposition~3.10 in \cite{BMS1}):
\begin{equation}
\lct(\fra)=\max_{b\in Q}\sum_ib_i.
\end{equation}
Note that since for every $j\leq r$ we have $\alpha_{i,j}>0$ for some $i$,
$Q$ is bounded.

It is more complicated to compute $\nu_f(p^e)$ for $e\geq 2$.
Let $c=\lct_0(f)$.
We show now that $\nu_f(p^e)=c(p^e-1)$ 
for all $e$ when $p\equiv 1$ (mod $N$),
where $N$ will be suitably chosen.

If $\lct_0(f)=1$, choose any $v\in Q\cap\QQ^r$ such that $\sum_iv_i=1$,
and if $\lct_0(f)<1$, let $v$ be one of the vertices of $Q$, such that
$\sum_iv_i=\max_{b\in Q}\sum_ib_i$. We take $N$ such that $Nv_i$
is an integer for all $i$. Note that we have $c=\sum_iv_i$.

If $p\equiv 1$ (mod $N$), then $a_i:=(p^e-1)v_i\in\NN$, and in order to
show that $\nu_f(p^e)\geq c(p^e-1)$, it is enough
to show
that $p$ does not divide $((p^e-1)c)!/a_1!\ldots a_r!$.
Therefore it is enough to check that
$$\lfloor c(p^e-1)/p^{e'}\rfloor=\sum_{i=1}^r\lfloor (p^e-1)v_i/p^{e'}\rfloor$$
for $1\leq e'\leq e$, where $\lfloor x\rfloor$ is the integral part of $x$.
This follows by an easy computation.

Dividing by $p^e$ and passing to the limit,
we deduce  that if $p\equiv 1$ (mod $N$),
then $\pt(f_p)\geq\lct_0(f)$. On the other hand,
 the reverse inequality
holds by Theorem~\ref{thm1'} (here this follows also from
the fact that $f\in \fra$ and 
$\lct_0(f)=\lct(\fra)$). Note that using Proposition~\ref{prop3},
we deduce that $\nu_f(p^e)=c(p^e-1)$ for all $e$ and all $p$ as above.
In particular,  we see that $\pt(f_p)$ exhibits the
behavior described in Problem~\ref{cong1}.
\end{example}

\begin{example}
Let $f=x^2+y^3$, so $\lct_0(f)=\frac{1}{2}+\frac{1}{3}=\frac{5}{6}$.
Moreover, $b_{f,0}$ has simple roots $-\frac{5}{6},-1,-\frac{7}{6}$
(see \cite{kashiwara1}, Example~6.19).
We give below the list of F-pure thresholds of $f_p$. We see that
for $p\neq 2$, $3$, the behavior depends on the congruence class of $p$
mod $3$. Note that we see the behavior described in Problems~\ref{cong1},
\ref{limit_p} and \ref{classes}
(when $J=\frmm$). Moreover, we obtain by our procedure all the roots of 
$b_{f,0}(s)$.
\begin{enumerate}
\item If $p=2$, then $\pt(f_p)=\frac{1}{2}$.

\item If $p=3$, we have $\pt(f_p)=\frac{2}{3}$.

\item If If $p\equiv 1$ (mod $3$), then $\pt(f_p)=\frac{5}{6}$.
We have $\nu_f(p^e)=\frac{5}{6}(p^e-1)$ for all $e$,
so this gives the root $-\frac{5}{6}$ of $b_{f,0}(s)$.

\item If $p\equiv 2$ (mod $3$) and $p\neq 2$, then
$\pt(f_p)=\frac{5}{6}-\frac{1}{6p}$, so
$$\nu_f(p^e)=
\begin{cases}
\frac{5}{6}p-\frac{7}{6}&\text{if $e=1$,}\\
\frac{5}{6}p^e-\frac{1}{6}p^{e-1}-1&\text{if $e\geq 2$.}
\end{cases}$$
Therefore we get the roots $-\frac{7}{6}$ and $-1$ of $b_{f,0}(s)$.
\end{enumerate}
\end{example}

\begin{example}\label{example3}
Let $f=x^2+y^7$. The log canonical threshold is given by
$\lct_0(f)=\frac{1}{2}+\frac{1}{7}=\frac{9}{14}$. All the roots of
$b_{f,0}$ are simple, and they are 
$$-\frac{9}{14},-\frac{11}{14},-\frac{13}{14},-1,-\frac{15}{14},
-\frac{17}{14},-\frac{19}{14}$$
(see \cite{kashiwara1}, Example~6.19).

We assume that $p\neq 2$, $7$ and we give the description
of the F-pure thresholds and of the functions $\nu_f(p^e)$.
The behavior depends on the congruence class of $p$ mod $7$.
We see again the behavior described in Problems~\ref{cong1},
\ref{limit_p} and \ref{classes}. In addition, we get all the roots of $b_{f,0}$
by our procedure.

\begin{enumerate}
\item If $p\equiv 1$ (mod $7$), then $\pt(f_p)=\frac{9}{14}$,
and $\nu_f(p^e)=\frac{9}{14}(p^e-1)$ for all $e$. This gives the
root $-\frac{9}{14}$ of $b_{f,0}(s)$.

\item If $p\equiv 2$ (mod $7$), then $\pt(f_p)=\frac{9}{14}-\frac{1}{14p^2}$.
Hence
$$\nu_f(p^e)=
\begin{cases}
\frac{9}{14}p-\frac{11}{14} & \text{if $e=1$,}\\
\frac{9}{14}p^2-\frac{15}{14} & \text{if $e=2$,}\\
\frac{9}{14}p^e-\frac{1}{14}p^{e-2}-1 & \text{if $e\geq 3$.}
\end{cases}$$
This gives the roots $-\frac{11}{14}$, $-\frac{15}{14}$ and $-1$
of $b_{f,0}(s)$.

\item If $p\equiv 3$ (mod $7$), then $\pt(f_p)=\frac{9}{14}-\frac{5}{14p^3}$.
Therefore 
$$\nu_f(p^e)=
\begin{cases}
\frac{9}{14}p-\frac{13}{14} & \text{if $e=1$,}\\
\frac{9}{14}p^2-\frac{11}{14} & \text{if $e=2$,}\\
\frac{9}{14}p^3-\frac{19}{14} & \text{if $e=3$,}\\
\frac{9}{14}p^e-\frac{5}{14}p^{e-3}-1 &\text{if $e\geq 4$.}
\end{cases}$$
We get the roots $-\frac{13}{14}$, $-\frac{11}{14}$, $-\frac{19}{14}$
and $-1$ of $b_{f,0}(s)$.

\item If $p\equiv 4$ (mod $7$), then $\pt(f_p)=\frac{9}{14}-\frac{1}{14p}$.
Hence 
$$\nu_f(p^e)=
\begin{cases}
\frac{9}{14}p-\frac{15}{14} & \text{if $e=1$,}\\
\frac{9}{14}p^e-\frac{1}{14}p^{e-1}-1 & \text{if $e\geq 2$.}
\end{cases}$$
This gives the roots $-\frac{15}{14}$ and $-1$ of $b_{f,0}(s)$.

\item If $p\equiv 5$ (mod $7$), then $\pt(f_p)=\frac{9}{14}-\frac{3}{14p}$.
Therefore
$$\nu_f(p^e)=
\begin{cases}
\frac{9}{14}p-\frac{17}{14} & \text{if $e=1$,}\\
\frac{9}{14}p^e-\frac{3}{14}p^{e-1}-1 & \text{if $e\geq 2$.}
\end{cases}$$
This gives the roots $-\frac{17}{14}$ and $-1$ for $b_{f,0}(s)$.

\item If $p\equiv 6$ (mod $7$), then $\pt(f_p)=\frac{9}{14}-\frac{5}{14p}$.
Hence 
$$\nu_f(p^e)=
\begin{cases}
\frac{9}{14}p-\frac{19}{14} & \text{if $e=1$,}\\
\frac{9}{14}p^e-\frac{5}{14}p^{e-1}-1 & \text{if $e\geq 2$.}
\end{cases}$$
This gives the roots $-\frac{19}{14}$ and $-1$ of $b_{f,0}(s)$.
\end{enumerate}
\end{example}

\begin{example}\label{ex4}
Let $f=x^5+y^4+x^3y^2$. The following are the roots of $b_{f,0}(s)$
(see \cite{Yano})
$$-\frac{9}{20},-\frac{11}{20},-\frac{13}{20},-\frac{7}{10},
-\frac{17}{20},-\frac{9}{10},-\frac{19}{20},-1,-\frac{21}{20},
-\frac{11}{10},-\frac{23}{20},-\frac{13}{10},-\frac{27}{20}.$$
As in Example~\ref{nondegenerate}, since the exponents of the monomials
in $f$ satisfy that genericity condition, we can compute the
jumping coefficients of $f$ using Howald's description from \cite{howald}.
The ones in $(0,1]$ are
$$\frac{9}{20},\,\frac{13}{20},\,\frac{7}{10},\,\frac{17}{20},\,
\frac{9}{10},\,\frac{19}{20},\,1.$$
As pointed out in \cite{Saito}, the interest in this example comes from the
fact that there is
a root $\lambda\in (-1,0)$ of $b_{f,0}(s)$ 
such that $-\lambda$ is not a jumping coefficients of $f$ (namely 
$\lambda=-\frac{11}{20}$). Note also that $f$ has an isolated singularity
at the origin.

We show that we can get all roots of $b_{f,0}(s)$ by the procedure
described in Remark~\ref{roots}. Note however that it is not enough to 
consider only $\pt^{J}_f(p^e)$ for $J=\frmm$. We will use the notation in
Proposition~\ref{notation} to index the functions $\nu^J(-)$. In particular,
$\nu_1(p^e)=\nu_{f}(p^e)$.

We assume that $p\neq 2$, $5$ and we compute $\nu_1(p^e)$, depending on the
congruence class of $p$ mod $20$.

\begin{enumerate}
\item If $p\equiv 1$ (mod $20$), then $\nu_1(p^e)=\frac{9}{20}(p^e-1)$
for all $e\geq 1$.

\item If $p\equiv 3$ (mod $20$), then $\nu_1(p^e)=
\begin{cases}
\frac{9}{20}p-\frac{27}{20} & \text{if $e=1$,}\\
\frac{9}{20}p^e-\frac{7}{20}p^{e-1}-1 & \text{if $e\geq 2$.}
\end{cases}$

\item If $p\equiv 7$ (mod $20$), then $\nu_1(p^e)=
\begin{cases}
\frac{9}{20}p-\frac{23}{20} & \text{if $e=1$,}\\
\frac{9}{20}p^e-\frac{3}{20}p^{e-1}-1 & \text{if $e\geq 2$.}
\end{cases}$

\item If $p\equiv 9$ (mod $20$), then $\nu_1(p^e)=
\begin{cases}
\frac{9}{20}-\frac{21}{20} & \text{if $e=1$,}\\
\frac{9}{20}p^e-\frac{1}{20}p^{e-1}-1 &\text{if $e\geq 2$.}
\end{cases}$

\item If $p\equiv 11$ (mod $20$), then 
$$\nu_1(p^e)=
\begin{cases}
\left(\frac{9}{20}p-\frac{19}{20}\right)\cdot\frac{p^{e+1}-1}{p^2-1}
+\left(\frac{19}{20}p-\frac{9}{20}\right)\cdot\frac{p^e-p}{p^2-1}
& \text{if $e$ is odd,}\\
\frac{9}{20}(p^e-1)
& \text{if $e$ is even.}
\end{cases}$$

\item If $p\equiv 13$ (mod $20$), then $\nu_1(p^e)=
\begin{cases}
\frac{9}{20}p-\frac{17}{20} & \text{if $e=1$,}\\
\frac{9}{20}p^2-\frac{21}{20} & \text{if $e=2$,}\\
\frac{9}{20}p^e-\frac{1}{20}p^{e-2}-1 & \text{if $e\geq 3$.}
\end{cases}$

\item If $p\equiv 17$ (mod $20$), then $\nu_1(p^e)=
\begin{cases}
\frac{9}{20}p-\frac{13}{20} & \text{if $e=1$,}\\
\frac{9}{20}p^2-\frac{21}{20} & \text{if $e=2$,}\\
\frac{9}{20}p^e-\frac{1}{20}p^{e-2}-1 & \text{if $e\geq 3$.}
\end{cases}$

\smallskip

\item If $p\equiv 19$ (mod $20$), then $\nu_1(p^e)=
\left(\frac{9}{20}p-\frac{11}{20}\right)(1+p+\ldots+p^{e-1})$
for all $e\geq 1$.
\end{enumerate}

We see that in this way we have accounted for the 
following roots of $b_{f,0}(s)$:  
$-\frac{9}{20}$, $-\frac{11}{20}$, $-\frac{13}{20}$,
$-\frac{17}{20}$, $-\frac{19}{20}$, $-1$,
$-\frac{21}{20}$, $-\frac{23}{20}$, $-\frac{27}{20}$.
In order to get the other four roots, we need to compute also
some values of $\nu_3(p)$.

\begin{enumerate}
\item If $p\equiv 1$ (mod $10$), then $\nu_3(p)=\frac{7}{10}(p-1)$.

\smallskip

\item If $p\equiv 3$ (mod $10$), then $\nu_3(p)=\frac{7}{10}p-\frac{11}{10}$.

\smallskip

\item If $p\equiv 5$ (mod $10$), then $\nu_3(p)=\frac{7}{10}p-\frac{9}{10}$.

\smallskip

\item If $p\equiv 7$ (mod $10$), then $\nu_3(p)=\frac{7}{10}p-\frac{13}{10}$.
\end{enumerate}
Therefore we recover also the roots $-\frac{7}{10}$, $-\frac{9}{10}$,
$-\frac{11}{10}$ and $-\frac{13}{10}$ of $b_{f,0}(s)$.

Note that if $p\equiv 19$ (mod $20$), then
$$\pt(f_p)=\frac{9p-11}{20(p-1)},$$
so this gives an example when
 $\pt(f_p)$ is not a polynomial in $\frac{1}{p}$. 
\end{example}

\begin{example}\label{elliptic}
Let $f$ be a homogeneous polynomial of degree $n$
in $\ZZ[X_1,\ldots,X_n]$, defining a smooth hypersurface
$Y$ in ${\mathbb P}^{n-1}$. It is well-known that in this case
$\lct_0(f)=1$ (see \cite{kollar}). 
On the other hand, it follows from Proposition~\ref{one_value}
that $\pt(f_p)=1$ if and only if the action of the Frobenius morphism
on $H_{\frmm}^{n-1}(R/f_p)$ is injective.
Here $R={\mathbb F}_p[X_1,\ldots,X_n]_{(X_1,\ldots,X_n)}$.

This action is injective if and only if it is injective on the
socle of $H^{n-1}_{\frmm}(R/(f_p))$. We assume that $p$ is large enough,
so $Y_p$, the reduction mod $p$ of $Y$, is smooth. It follows that 
$\pt(f_p)=1$ if and only if 
the action induced by the Frobenius morphism on 
$H^{n-2}(Y_p,{\mathcal O}_{Y_p})$ is injective. 

If $n=3$, then $Y$ is an elliptic curve. In this case we see that
$\pt(f_p)=1$ if and only if $Y$ is not supersingular. There are two cases:
suppose first that $Y$ has complex multiplication (over $\CC$). 
In this case, $Y_p$
is supersingular if and only if $p$ is inert in the imaginary quadratic
CM field.

On the other hand, if $Y$ has no complex multiplication then Serre
\cite{Se} proved that the set of primes $p$ for which $Y_p$ is supersingular
has natural density zero in the set of all primes. This clearly suggests that 
the behaviour of $\pt(f_p)$ does not depend on the congruence of $p$
modulo some $N$, as in Problems~\ref{cong1} and \ref{limit_p}.
It would be interesting to compute explicitly the F-pure thresholds
at the primes where the curve is supersingular.

We mention a result of Elkies \cite{El} which is relevant in this
setting: it says that for every elliptic curve
 $Y$ as above,
there are infinitely many primes $p$ for which $Y_p$ is supersingular,
and therefore $\pt(f_p)\neq 1$.

As B.~Conrad and N.~Katz pointed out to us, if $K$ is the field
obtained by adjoining to $\QQ$ all points of order $\ell$ of $Y$
(for some odd prime $\ell$), then for every odd prime $p$
such that $p$ splits completely in $K$, the curve $Y_p$
is not supersingular. This provides an affirmative answer to 
Question~\ref{number_field} in this example. 
\end{example}

\medskip

\subsection*{Acknowledgements}
We are grateful to Johan de Jong, Lawrence Ein and Martin Olsson
for useful discussions. We are particularly indebted to Karen Smith for
suggesting to us Example~\ref{elliptic}, and to
Brian Conrad and Nick Katz for answering our questions on 
supersingular elliptic curves. 

While working on this project
the first author
was a Clay Mathematics Institute Research Fellow. 
Our work started during
the first author's visit to University of Tokyo. He is grateful 
to his host Yujiro Kawamata for his wonderful hospitality.
Part of this work was done during the second's author's stay at 
University of Michigan. He would like to express his deep gratitude
to Melvin Hochster for his hospitality and support.

\providecommand{\bysame}{\leavevmode \hbox \o3em

{\hrulefill}\thinspace}

\end{document}